# EXACT ESTIMATES FOR MOMENTS OF RANDOM BILINEAR FORMS

By R. Ibragimov[1], Sh. Sharakhmetov[2] and A. Cecen[3]

**Running head:** Exact estimates for moments of random bilinear forms

**Abstract.** The present paper concentrates on the analogues of Rosenthal's inequalities for ordinary and decoupled bilinear forms in symmetric random variables. More specifically, we prove the exact moment inequalities for these objects in terms of moments of their individual components. As a corollary of these results we obtain the explicit expressions for the best constant in the analogues of Rosenthal's inequality for ordinary and decoupled bilinear forms in identically distributed symmetric random variables in the case of the fixed number of random variables.

**Mailing address:** Department of Economics, Yale University, 28 Hillhouse Ave., New Haven, CT 06511; Phone: (203) 7772750; E-Mail: rustam.ibragimov@yale.edu

---

[1] Department of Economics, Yale University, 28 Hillhouse Ave., New Haven, CT 06511. E-Mail: rustam.ibragimov@yale.edu. Supported in part by a grant from Central Michigan University Graduate Studies and by a grant from the American Council of Teachers of Russian/American Council for Collaboration in Education and Language Study, with funds provided by the United States Information Agency

[2] Department of Probability Theory, Tashkent State Economics University, Tashkent, Uzbekistan

[3] Department of Economics, Central Michigan University, Mt. Pleasant, Michigan, USA

Supported in part by a grant from Central Michigan University Department of Economics

**1. Introduction.** In recent years, several studies have focused on moment and probability inequalities for multilinear forms and symmetric statistics (see, in particular, Serfling (1980), Krakowiak and Szulga (1986), McConell and Taqqu (1986), de la Pena (1992), de la Pena and Klass (1994), Koroljuk and Borovskikh (1994), de la Pena and Montgomery-Smith (1995), Sharakhmetov (1995, 1997), Ibragimov and Sharakhmetov (1996a, 1999), Borovskikh and Korolyuk (1997), Ibragimov (1997), Klass and Nowicki (1997a, b) and Gine et. al. (2000)). Interest in such inequalities is motivated by their applications in limit theorems, multiple stochastic integration, harmonic analysis, operator theory, quantum mechanics, theory of income inequality and species' diversity measurement, etc. (see, in addition to the above-mentioned papers, Bonami (1970), Rosinski and Szulga (1982), Sjorgen (1982), Rosinski and Woyczynski (1984, 1986), Cambanis et al. (1985) and Kwapien and Woyczinski (1992)). Furthermore, the bounds on moments for symmetric statistics can also be applied in investment theory and in testing for chaos in time series data based on the notion of correlation integral, which has the form of symmetric statistics (see Cecen and Erkal (1996a, b)).

In the case of linear statistics (sums of independent random variables (r.v.'s)), the exact moment estimates are given by the well-known Khintchine, Marcinkiewicz-Zygmund and Rosenthal inequalities (see Khintchine (1923), Marcinkiewicz and Zygmund (1937), Rosenthal (1970)). Let us remind the latter ones ($A_i(\cdot)$, $B_i(\cdot)$ denote constants depending on parameters in parentheses only).

**Theorem 1.** If $\xi_1,\ldots,\xi_n$ are independent mean zero r.v.'s with finite $t$th moment, $2<t<\infty$, then

$$A_1(t)max\left(\sum_{i=1}^{n}E|\xi_i|^t,\left(\sum_{i=1}^{n}E\xi_i^2\right)^{t/2}\right)\le E\left|\sum_{i=1}^{n}\xi_i\right|^t\le B_1(t)max\left(\sum_{i=1}^{n}E|\xi_i|^t,\left(\sum_{i=1}^{n}E\xi_i^2\right)^{t/2}\right). \qquad (1)$$

The exact upper constants in inequality (1) (case $t=2m$) and in its analogue for nonnegative r.v.'s were found in Ibragimov and Sharakhmetov (1996b, 1998a, b). The best constant in inequality (1) for symmetric r.v.'s was independently found by Figiel et al. (1997) and Ibragimov and Sharakhmetov (1995, 1997). The results obtained by Ibragimov and Sharakhmetov (1996b, 1997, 1998a, b) and their proofs were presented in Ibragimov (1997). Concerning refinements and extensions of Rosenthal's inequalities and related problems see also Prokhorov (1962), Nagaev and Pinelis (1977), Pinelis (1980, 1994), Pinelis and Utev (1984), Johnson et al. (1985), Utev (1985), Talagrand (1989), Hitczenko (1990, 1994), Nagaev (1990, 1998), Kwapien and Szulga (1991) and Peshkir and Shiryaev (1995).

Recently, Sharakhmetov (1995, 1997), Ibragimov and Sharakhmetov (1996a, 1998a, 1999, 2000) (see also Ibragimov (1997)), Klass and Nowicki (1997a, b) and Gine et. al. (2000) obtained analogues of Rosenthal's inequality (1) and its analogue for nonnegative r.v.'s in the case of symmetric statistics. Ibragimov and Sharakhmetov (2000) also showed the significance of each term in the analogues of Rosenthal's bounds for symmetric statistics of arbitrary order. Ibragimov (1997) showed that the best constants in the analogues of Rosenthal's inequalities grow not slower than $(t/ln\,t)^m$, as $t\to\infty$, where $m$ is the order of a symmetric statistic. Gine et. al. (2000) showed that the actual rate of growth of the above constants is $(t/\ln t)^m$.

The qualitative difference of the results on Rosenthal's inequalities for nonlinear statistics from the linear case is the exact constants in them are unknown yet. The main goal of the present paper is to fill partially this gap in the case of bilinear forms. More specifically, we obtain the explicit expressions for the best constant in the analogues of Rosenthal's inequalities for ordinary and decoupled bilinear forms in identically distributed symmetric r.v.'s in the case of fixed number of r.v.'s. The proof of the expressions for the best constants in the non-linear analogues of Rosenthal inequalities is based on a theorem, which extends the extremal results obtained in Utev (1985) and Ibragimov and Sharakhmetov (1996b, 1997) in the case of bilinear forms and gives the exact estimates for moments of random bilinear forms in terms of moment characteristics of their particular components. To our knowledge, this theorem and its proof are the first attempt to apply methods which were used to investigate the extremal problems in moment inequalities for sums of independent r.v.'s for non-linear statistics. The results obtained in the present paper can be extended to the case of nonnegative random variables, multilinear forms of arbitrary order and generalized moments; these extensions will be presented elsewhere.

**2. Main results.** Let $t>2$, $X_1, Y_1, X_2, Y_2, \ldots, X_n, Y_n$ be independent symmetric r.v.'s with finite $t$th moment. Let $a_i \geq 0$, $b_i \geq 0$, $c_i \geq 0$, $d_i \geq 0$, $a_i^t \leq b_i$, $c_i^t \leq d_i, i = 1, \ldots, n$. Set

$$(X, n) = (X_1, \dots, X_n),\ (Y, n) = (Y_1, \dots, Y_n)$$

$$M_1(n,a,b) = \{(X,n) : EX_i^2 = a_i^2, E|X_i|^t = b_i,\ i = 1,\dots,n\}$$

$$M_1(n,c,d) = \{(Y,n) : EY_i^2 = c_i^2, E|Y_i|^t = d_i,\ i = 1,\dots,n\}$$

$$M_2(n,a,b) = \{(X,n) : EX_i^2 \le a_i^2, E|X_i|^t \le b_i,\ i = 1,\dots,n\}$$

$$M_2(n,c,d) = \{(Y,n) : EY_i^2 \le c_i^2, E|Y_i|^t \le d_i,\ i = 1,\dots,n\}$$

Let $U_i(a_i, b_i, t),\ V_i(c_i, d_i, t),\ i = 1,\dots,n,$ be independent r.v.'s such that

$$P(U_i(a_i, b_i, t) = 0) = 1 - (a_i^t / b_i)^{2/(t-2)}$$

$$P(U_i(a_i, b_i, t) = \pm(b_i / a_i^2)^{1/(t-2)}) = (1/2)(a_i^t / b_i)^{2/(t-2)}$$

$$P(V_i(c_i, d_i, t) = 0) = 1 - (c_i^t / d_i)^{2/(t-2)}$$

$$P(V_i(c_i, d_i, t) = \pm(d_i / c_i^2)^{1/(t-2)}) = (1/2)(c_i^t / d_i)^{2/(t-2)}$$

and let $U_i$, $V_i$, $i=1,...,n$, be independent r.v.'s with distribution

$$P(U_i=\pm 1)=P(V_i=\pm 1)=1/2,\ i=1,...,n$$

The following theorem extends the results obtained in Utev (1985) and Ibragimov and Sharakhmetov (1997) on the non-linear case and gives the explicit bounds for moments of random bilinear forms in terms of moment characteristics of their particular components.

**Theorem 4.** If $2<t<4$, then

$$\sup_{(X,n)\in M_k(n,a,b)} E\left|\sum_{1\le i<j\le n} X_i X_j\right|^t = \sum_{1\le i<j\le n}(b_i-a_i^t)(b_j-a_j^t)+$$

$$+\sum_{i=1}^n (b_i-a_i^t)E\left|\sum_{j\ne i}^n a_j U_j\right|^t + E\left|\sum_{1\le i<j\le n} a_i a_j U_i U_j\right|^t \tag{2}$$

$$\sup_{\substack{(X,n)\in M_k(n,a,b),\\ (Y,n)\in M_l(n,c,d)}} E\left|\sum_{1\le i,j\le n,\, i\ne j} X_i Y_j\right|^t = \sum_{\substack{1\le i,j\le n\\ i\ne j}}(b_i-a_i^t)(d_j-c_j^t)+$$

$$+\sum_{j=1}^n (d_j-c_j^t)E\left|\sum_{\substack{i=1,\\ i\ne j}}^n a_i U_i\right|^t+\sum_{i=1}^n (b_i-a_i^t)E\left|\sum_{j=1,\ j\ne i}^n c_j V_j\right|^t+$$

$$+E\left|\sum_{1\le i,j\le n,\ i\ne j} a_i c_j U_i V_j\right|^t,\ \ k,l=1,2 \tag{3}$$

If $3 \le t < 4$, then

$$\inf_{(X,n)\in M_1(n,a,b)} E\left|\sum_{1\le i<j\le n} X_i X_j\right|^t = E\left|\sum_{1\le i<j\le n} U_i(a_i,b_i,t)U_j(a_j,b_j,t)\right|^t \tag{4}$$

$$\inf_{\substack{(X,n)\in M_1(n,a,b),\\ (Y,n)\in M_1(n,c,d)}} E\left|\sum_{\substack{1\le i,j\le n,\\ i\ne j}} X_i Y_j\right|^t = E\left|\sum_{\substack{1\le i,j\le n,\\ i\ne j}} U_i(a_i,b_i,t)V_j(c_j,d_j,t)\right|^t \tag{5}$$

If $t \ge 4$, then

$$\sup_{(X,n)\in M_k(n,a,b)} E\left|\sum_{1\le i<j\le n} X_i X_j\right|^t = E\left|\sum_{1\le i<j\le n} U_i(a_i,b_i,t)U_j(a_j,b_j,t)\right|^t \tag{6}$$

$$\sup_{\substack{(X,n)\in M_k(n,a,b),\\ (Y,n)\in M_l(n,c,d)}} E\left|\sum_{\substack{1\le i,j\le n,\\ i\ne j}} X_i Y_j\right|^t = E\left|\sum_{\substack{1\le i,j\le n,\\ i\ne j}} U_i(a_i,b_i,t)V_j(c_j,d_j,t)\right|^t ,k,l=1,2 \tag{7}$$

$$\inf_{(X,n)\in M_1(n,a,b)} E\left|\sum_{1\le i<j\le n} X_i X_j\right|^t = \sum_{1\le i<j\le n} (b_i - a_i^t)(b_j - a_j^t) +$$

$$+\sum_{i=1}^{n} (b_i - a_i^t) E\left|\sum_{j=1,\, j\ne i}^{n} a_j U_j\right| + E\left|\sum_{1\le i<j\le n} a_i a_j U_i U_j\right|^t \tag{8}$$

$$\inf_{\substack{(X,n)\in M_1(n,a,b),\\ (Y,n)\in M_1(n,c,d)}} E\left|\sum_{\substack{1\le i,j\le n,\\ i\ne j}} X_i Y_j\right|^t = \sum_{\substack{1\le i,j\le n,\\ i\ne j}} (b_i - a_i^t)(d_j - c_j^t) +$$

$$+\sum_{j=1}^{n} (d_j - c_j^t) E\left|\sum_{\substack{i=1,\\ i\ne j}}^{n} a_i U_i\right|^t + \sum_{i=1}^{n} (b_i - a_i^t) E\left|\sum_{\substack{j=1,\\ j\ne i}}^{n} c_j V_j\right|^t +$$

$$+E\left|\sum_{\substack{1\le i,j\le n,\\ i\ne j}} a_i c_j U_i V_j\right|^t \qquad (9)$$

<u>Remark.</u> The expressions in relations (2)-(9) are of a simple structure and their values can be easily calculated for given sequences $a_i$, $b_i$, $c_i$, $d_i$, $i = 1,...,n$.

Let us fix $t>2$ and $n \ge 1$. From the results obtained in Ibragimov and Sharakhmetov (1999) and decoupling theorems for symmetric statistics (see McConell and Taqqu (1986) and de la Pena and Montgomery-Smith (1995)) it follows that for all independent identically distributed symmetric r.v.'s $X_1,...,X_n$, $\overline{X}_1,...,\overline{X}_n$ with finite $t$th moment the following Rosenthal-type inequalities are true ($C_n^2 = n(n-1)/2$):

$$E\left|\sum_{1\le i<j\le n} X_i X_j\right|^t \le B_4(t,n)\ \max(C_n^2 (E|X_1|^t)^2, (C_n^2)^{t/2} (EX_1^2)^t) \qquad (10)$$

$$E\left|\sum_{1\le i<j\le n} X_i X_j\right|^t \le B_5(t,n)\ \max(n^2 (E|X_1|^t)^2, n^t (EX_1^2)^t) \qquad (11)$$

$$E\left|\sum_{1\le i\ne j\le n} X_i \overline{X}_j\right|^t \le B_6(t,n) max(C_n^2(E|X_1|^t)^2, (C_n^2)^{t/2}(EX_1^2)^t) \tag{12}$$

$$E\left|\sum_{1\le i\ne j\le n} X_i \overline{X}_j\right|^t \le B_7(t,n) max(n^2(E|X_1|^t)^2, n^t(EX_1^2)^t) \tag{13}$$

The following theorems give the explicit expressions for the best constants in inequalities (10) and (11).

**Theorem 5.** The exact constant in inequality (10) is given by

$$B_4^*(t,n) = C_n^2(1/(C_n^2)^{1/2} - 1/(C_n^2)^{t/2})^2 +$$

$$+(1/(C_n^2)^{1/2} - 1/(C_n^2)^{t/2})n/(C_n^2)^{t/4} E\left|\sum_{i=2}^n U_i\right|^t +$$

$$+E\left|\sum_{1\le i<j\le n} U_i U_j /(C_n^2)^{1/2}\right|^t, \ 2 < t < 4 \tag{14}$$

$$B_4^*(t,n) = E\left|\sum_{1\le i<j\le n} U_i(1/(C_n^2)^{1/4}, 1/(C_n^2)^{1/2}, t)\right.$$

$$\left.\times U_j(1/(C_n^2)^{1/4}, 1/(C_n^2)^{1/2}, t)\right|^t, \ t \ge 4 \tag{15}$$

**Theorem 6.** The exact constant in inequality (11) is given by

$$B_5^*(t,n) = C_n^2(1/n - 1/n^t)^2 + (1/n^{t/2} - 1/n^{3t/2-1})E\left|\sum_{i=2}^{n} U_i\right|^t +$$

$$+ E\left|\sum_{1 \le i < j \le n} U_i U_j / n\right|^t, \ 2 < t < 4 \qquad (16)$$

$$B_5^*(t,n) = E\left|\sum_{1 \le i < j \le n} U_i(1/n^{1/2},\ 1/n,\ t)U_j(1/n^{1/2},\ 1/n,\ t)\right|^t, \ \ t \ge 4 \qquad (17)$$

Theorems 7 and 8 below give the explicit expressions for the exact constants in inequalities (12) and (13).

**Theorem 7.** The exact constant in inequality (12) is given by

$$B_6^*(t,n) = 2C_n^2(1/(C_n^2)^{1/2} - 1/(C_n^2)^{t/2})^2 +$$

$$+2(1/(C_n^2)^{1/2} - 1/(C_n^2)^{t/2})n/(C_n^2)^{t/4}E\left|\sum_{i=2}^{n} U_i\right|^t +$$

$$+ E\left|\sum_{1 \le i, j \le n,\ i \ne j} U_i V_j /(C_n^2)^{1/2}\right|^t, \ 2 < t < 4 \qquad (18)$$

$$B_6^*(t,n) =$$

$$E\left|\sum_{1\le i\ne j\le n} U_i(1/(C_n^2)^{1/4}, 1/(C_n^2)^{1/2}, t)\right.$$

$$\left.\times V_j(1/(C_n^2)^{1/4}, 1/(C_n^2)^{1/2}, t)\right|^t, \quad t \ge 4 \tag{19}$$

**Theorem 8.** The exact constant in inequality (13) is given by

$$B_7^*(t,n) = 2C_n^2(1/n - 1/n^t)^2 + 2(1/n^{t/2} - 1/n^{3t/2-1})E\left|\sum_{i=2}^n U_i\right|^t +$$

$$+E\left|\sum_{\substack{1\le i,j\le n,\\ i\ne j}} U_i V_j / n\right|^t, \quad 2 < t < 4 \tag{20}$$

$$B_7^*(t,n) = E\left|\sum_{1\le i,j\le n,\, i\ne j} U_i(1/n^{1/2}, 1/n, t) V_j(1/n^{1/2}, 1/n, t)\right|^t,$$

$$t \ge 4 \tag{21}$$

**3. Preliminaries.** Let us formulate some auxiliary steps needed for the proof of the theorems.

**Lemma 1.** If $2 < t < 4$, $z_1 \ge 0, z_2 \in \mathbf{R}$, $a \ge 0$, $b \ge 0$, $a^t \le b$ , $X$ is a symmetric r.v. with $EX^2 \le a^2$, $E|X|^t \le b$, then

$$E|z_1 X + z_2|^t - bz_1^t \le E|az_1 U + z_2|^t - a^t z_1^t \tag{22}$$

Proof. It suffices to consider the case $z_1 \neq 0$. From Lemma 5 in Ibragimov and Sharakhmetov (1997) it follows that

$$E\left|X + z_2 / z_1\right|^t - b \leq E\left|aU + z_2 / z_1\right|^t - a^t \qquad (23)$$

Multiplying (23) by $z_1^t$ we obtain (22).

Applying Lemma 7 in Ibragimov and Sharakhmetov (1997) and Lemmas 7.3 and 7.4 in Utev (1985) analogously to the proof of Lemma 1 above we easily obtain the following Lemmas 2-4.

**Lemma 2.** If $3 \leq t < 4$, $z_1, z_2 \in \mathbf{R}$, $a \geq 0$, $b \geq 0$, $a^t \leq b$, $X$ is a symmetric r.v. with $EX^2 = a^2$, $E|X|^t = b$, then

$$E\left|z_1 X + z_2\right|^t \geq E\left|z_1 U(a,b,t) + z_2\right|^t$$

**Lemma 3.** If $t \geq 4$, $z_1, z_2 \in \mathbf{R}$, $a \geq 0$, $b \geq 0$, $a^t \leq b$, $X$ is a symmetric r.v. with $EX^2 \leq a^2$, $E|X|^t \leq b$, then

$$E\left|z_1 X + z_2\right|^t \leq E\left|z_1 U(a,b,t) + z_2\right|^t$$

**Lemma 4.** If $t \ge 4$, $z_1 \ge 0, z_2 \in \mathbf{R}$, $a \ge 0$, $b \ge 0$, $a^t \le b$, $X$ is a symmetric r.v. with $EX^2 = a^2$, $E|X|^t = b$, then

$$E\left|z_1 X + z_2\right|^t - bz_1^t \ge E\left|az_1 U + z_2\right|^t - a^t z_1^t.$$

**Lemma 5.** Let $1 \le k \le n$, $X_1,\ldots, X_{k-1}, U_k, X_{k+1},\ldots, X_n$ be independent r.v.'s with $E\left|X_i\right|^t < \infty$, $i = 1,\ldots, n$, $i \ne k$, $a_k, b_k \ge 0$, $a_k^t \le b_k$, $c_i \in \mathbf{R}$, $i = 1,\ldots, k-1$, and let $F_1$ be the set of symmetric r.v.'s $X_k$ independent of $X_1,\ldots, X_{k-1}, X_{k+1},\ldots, X_n$ and satisfying the conditions $EX_k^2 \le a_k^2$, $E\left|X_k\right|^t \le b_k$, $F_2$ be the subset of $F_1$ consisting of r.v.'s $X_k$ such that $EX_k^2 = a_k^2$, $E\left|X_k\right|^t = b_k$. If $2 < t < 4$, then

$$\sup_{X_k \in F_l} \Big( \sum_{i=1}^{k-1} c_i E\Big|\sum_{\substack{j=1 \\ j \ne i}}^{n} X_j\Big|^t + E\Big|\sum_{1 \le i < j \le n} X_i X_j\Big|^t \Big)$$

$$= \sum_{i=1}^{k-1} c_i E\Big|a_k U_k + \sum_{\substack{j=1 \\ j \ne i,k}}^{n} X_j\Big|^t + \sum_{i=1}^{k-1} c_i (b_k - a_k^t) +$$

$$+(b_k - a_k^t)E\left|\sum_{\substack{j=1\\ j\neq k}}^{n} X_j\right|^t + E\left|a_k U_k(\sum_{\substack{j=1\\ j\neq k}}^{n} X_j) + \sum_{\substack{1\le i<j\le n\\ i,j\neq k}} X_i X_j\right|^t, \; l=1,2$$

If $t \geq 4$, then

$$\inf_{X_k\in F_2}\left(\sum_{i=1}^{k-1} c_i E\left|\sum_{\substack{j=1\\ j\neq i}}^{n} X_j\right|^t + E\left|\sum_{1\le i<j\le n} X_i X_j\right|^t\right)$$

$$=\sum_{i=1}^{k-1} c_i E\left|a_k U_k + \sum_{\substack{j=1\\ j\neq i,k}}^{n} X_j\right|^t + \sum_{i=1}^{k-1} c_i(b_k - a_k^t) +$$

$$+(b_k - a_k^t)E\left|\sum_{\substack{j=1\\ j\neq k}}^{n} X_j\right|^t + E\left|a_k U_k(\sum_{\substack{j=1\\ j\neq k}}^{n} X_j) + \sum_{\substack{1\le i<j\le n\\ i,j\neq k}} X_i X_j\right|^t$$

Proof. From Lemmas 1 and 4 above and Lemma 5 in Ibragimov and Sharakhmetov (1997) it follows that it suffices to find a sequence of r.v.'s $X_{mk}$, $m=1,2,\ldots,$ independent of $X_1,\ldots,X_{k-1},X_{k+1},\ldots,X_n$ and satisfying the conditions

$EX_{mk}^2 = a_k^2,\; E|X_{mk}|^t = b_k,$

$$\lim_{m\to\infty} E\left|X_{mk} + \sum_{\substack{j=1\\ j\neq i,k}}^{n} X_j\right|^t = E\left|a_k U_k + \sum_{\substack{j=1\\ j\neq i,k}}^{n} X_j\right|^t + b_k - a_k^t,\; i=1,\ldots,k-1 \qquad (24)$$

$$\lim_{m\to\infty} E\left|X_{mk}(\sum_{\substack{j=1\\ j\neq k}}^{n} X_j) + \sum_{\substack{1\le i<j\le n\\ i,j\neq k}} X_i X_j\right|^t = (b_k - a_k^t)E\left|\sum_{\substack{j=1\\ j\neq k}}^{n} X_j\right|^t +$$

$$+E\left|a_k U_k(\sum_{\substack{j=1\\ j\neq k}}^{n} X_j) + \sum_{\substack{1\le i<j\le n\\ i,j\neq k}} X_i X_j\right|^t \qquad (25)$$

If $b_k = a_k^t$, then one can take $X_{mk} = a_k U_k$. Let $a_k^t < b_k$. Set $\delta_m = 1/m$, $P(X_{mk} = \pm a_k) = 1/2(1-\delta_m)$, $P(X_{mk} = \pm b_{mk}) = 1/2\delta^*_{mk}$, $\delta^*_{mk} = a_k^2 \delta_m / b_{mk}^2$, $P(X_{mk} = 0) = \delta_m - \delta^*_{mk}$, $b_{mk} = ((b_k - a_k^t(1-\delta_m))/a_k^2\delta_m)^{1/(t-2)}$, $m = 1,2,...$

Then

$$b_{mk} \ge a_k,\ 0 \le \delta^*_{mk} \le \delta_m,$$

$$EX^2_{mk} = a_k^2,\ E\left|X_{mk}\right|^t = b_k,\ m = 1,2,... \qquad (26)$$

$$\delta_m \to 0,\ b_{mk} \to \infty,\ b^t_{mk}\delta^*_{mk} \to b_k - a_k^t,\ m \to \infty$$

From (26) and the proof of Lemma 7.6 in Utev (1985) it follows that relations (24) are valid.

Let us prove that (25) is true. We have

$$E\left|X_{mk}(\sum_{\substack{j=1,\\ j\neq k}}^{n} X_j)+\sum_{\substack{1\le i<j\le n,\\ i,j\neq k}} X_iX_j\right|^t = E\left|a_kU_k(\sum_{\substack{j=1,\\ j\neq k}}^{n} X_j)+\sum_{\substack{1\le i<j\le n,\\ i,j\neq k}} X_iX_j\right|^t(1-\delta_m)+$$

$$+E\left|\sum_{\substack{1\le i<j\le n,\\ i,j\neq k}} X_iX_j\right|^t(\delta_m-\delta^*_{mk})+(E\left|b_{mk}U_k(\sum_{\substack{j=1,\\ j\neq k}}^{n} X_j)+\sum_{\substack{1\le i<j\le n,\\ i,j\neq k}} X_iX_j\right|^t -$$

$$-b^t_{mk}E\left|\sum_{\substack{j=1\\ j\neq k}}^{n} X_j\right|^t)\delta^*_{mk}+b^t_{mk}\delta^*_{mk}E\left|\sum_{\substack{j=1,\\ j\neq k}}^{n} X_j\right|^t$$

From (26) it follows that for the proof of (25) it suffices to check that

$$(E\left|b_{mk}U_k(\sum_{\substack{j=1,\\ j\neq k}}^{n} X_j)+\sum_{\substack{1\le i<j\le n,\\ i,j\neq k}} X_iX_j\right|^t - b^t_{mk}E\left|\sum_{\substack{j=1\\ j\neq k}}^{n} X_j\right|^t)\delta^*_{mk}\to 0,\ m\to\infty$$

This follows from the fact that $b^t_{mk}\delta^*_{mk}$ converges and that, on the strength of the inequality $\left||x+y|^t-|x|^t\right|\le 2^t t(|x|^{t-1}|y|+|y|^t)$, $x,y\in\mathbf{R}$, $t\ge 1$ (see Lemma 7.5. in Utev (1984)), and the dominated convergence principle,

$$\lim_{m\to\infty} E\left|U_k(\sum_{\substack{j=1,\\ j\neq k}}^{n} X_j) + \sum_{\substack{1\le i<j\le n,\\ i,j\neq k}} X_i X_j / b_{mk}\right|^t = E\left|U_k(\sum_{\substack{j=1,\\ j\neq k}}^{n} X_j)\right|^t = E\left|\sum_{\substack{j=1,\\ j\neq k}}^{n} X_j\right|^t$$

Arguing analogously with the proof of Lemma 5, we easily obtain the following

**Lemma 6.** Let $1 \le k \le n$, $X_1,..., X_{k-1}, U_k, X_{k+1},..., X_n, Y_1,..., Y_n$ be independent r.v.'s with $E|X_i|^t < \infty$, $i=1,...,n,\ i\neq k$, $E|Y_i|^t < \infty$, $i=1,...,n,\ a_k, b_k \ge 0$, $a_k^t \le b_k$, $c_i \in \mathbf{R}$, $i=1,...,k-1$, and let $G_1$ be the set of symmetric r.v.'s $X_k$ independent of $X_1,..., X_{k-1}, X_{k+1},..., X_n, Y_1,..., Y_n$ and satisfying the conditions $EX_k^2 \le a_k^2$, $E|X_k|^t \le b_k$, $G_2$ be the subset of $G_1$ consisting of r.v.'s $X_k$ such that $EX_k^2 = a_k^2$, $E|X_k|^t = b_k$. If $2 < t < 4$, then

$$\sup_{X_k\in G_l} (\sum_{i=1}^{n} c_i E\left|\sum_{j=1,\ j\neq i}^{n} X_j\right|^t + E\left|\sum_{1\le i\neq j\le n} X_i Y_j\right|^t) =$$

$$= \sum_{i=1}^{n} c_i E\left|a_k U_k + \sum_{j=1,\ j\neq i,k}^{n} X_j\right|^t + \sum_{i=1}^{n} c_i (b_k - a_k^t) +$$

$$+ (b_k - a_k^t) E\left|\sum_{j=1,\ j\neq k}^{n} Y_j\right|^t + E\left|a_k U_k (\sum_{j=1,\ j\neq k}^{n} Y_j) + \sum_{1\le i\neq j\le n,\ i\neq k} X_i Y_j\right|^t,\quad l=1,2$$

If $t \ge 4$, then

$$\inf_{X_k \in G_2} \Big( \sum_{i=1}^{n} c_i E \Big| \sum_{j=1,\ j\neq i}^{n} X_j \Big|^t + E \Big| \sum_{1\leq i\neq j\leq n} X_i Y_j \Big|^t \Big) =$$

$$= \sum_{i=1}^{n} c_i E \Big| a_k U_k + \sum_{j=1,\ j\neq i,k}^{n} X_j \Big|^t + \sum_{i=1}^{n} c_i (b_k - a_k^t) +$$

$$+ (b_k - a_k^t) E \Big| \sum_{\substack{j=1 \\ j\neq k}}^{n} Y_j \Big|^t + E \Big| a_k U_k \Big( \sum_{j=1,\ j\neq k}^{n} Y_j \Big) + \sum_{1\leq i\neq j\leq n,\ i,j\neq k} X_i Y_j \Big|^t$$

**4. Proofs of the theorems.**

**Proof of theorem 3.** Relations (4)-(7) easily follow from Lemmas 2 and 3 by induction. Let us prove (2). Let $2 < t < 4$, $1 \leq k \leq n$, $U_1,..., U_{k-1}$, $X_{k+1},..., X_n$ be independent symmetric r.v.'s, $E|X_i|^t < \infty$, $i = k+1,...,n$, $a_i \geq 0$, $b_i \geq 0$, $a_i^t \leq b_i$, $i = 1,...,k$. Denote by $H_1$ the set of symmetric r.v.'s $X_k$ independent of $U_1,..., U_{k-1}$, $X_{k+1},..., X_n$ and satisfying the conditions $EX_k^2 \leq a_k^2$, $E|X_k|^t \leq b_k$, and by $H_2$ the subset of $H_1$ consisting of r.v.'s $X_k$ such that $EX_k^2 = a_k^2$, $E|X_k|^t = b_k$. On the strength of Lemma 5 we have

$$\sup_{X_k \in H_l} \Big( \sum_{1\le i<j\le k-1} (b_i - a_i^t)(b_j - a_j^t) + \sum_{i=1}^{k-1} (b_i - a_i^t) E\Big|\sum_{\substack{j=1\\ j\ne i}}^{k-1} a_j U_j + \sum_{j=k}^{n} X_j\Big|^t +$$

$$+E\Big|\sum_{i=1}^{k-1} a_i U_i (\sum_{j=i+1}^{k-1} a_j U_j + \sum_{j=k}^{n} X_j) + \sum_{i=k}^{n-1} X_i (\sum_{j=i+1}^{n} X_j)\Big|^t \Big) =$$

$$= \sum_{1\le i<j\le k-1} (b_i - a_i^t)(b_j - a_j^t) + \sum_{i=1}^{k-1} (b_i - a_i^t) E\Big|\sum_{\substack{j=1\\ j\ne i}}^{k} a_j U_j + \sum_{j=k+1}^{n} X_j\Big|^t +$$

$$+\sum_{i=1}^{k-1} (b_i - a_i^t)(b_k - a_k^t) + (b_k - a_k^t) E\Big|\sum_{j=1}^{k-1} a_j U_j + \sum_{j=k+1}^{n} X_j\Big|^t +$$

$$+E\Big|\sum_{i=1}^{k} a_i U_i (\sum_{j=i+1}^{k} a_j U_j + \sum_{j=k+1}^{n} X_j) + \sum_{i=k+1}^{n-1} X_i (\sum_{j=i+1}^{n} X_j)\Big|^t \Big) =$$

$$= \sum_{1\le i<j\le k} (b_i - a_i^t)(b_j - a_j^t) + \sum_{i=1}^{k} (b_i - a_i^t) E\Big|\sum_{\substack{j=1\\ j\ne i}}^{k} a_j U_j + \sum_{j=k+1}^{n} X_j\Big|^t +$$

$$+E\Big|\sum_{i=1}^{k} a_i U_i (\sum_{j=i+1}^{k} a_j U_j + \sum_{j=k+1}^{n} X_j) + \sum_{i=k+1}^{n-1} X_i (\sum_{j=i+1}^{n} X_j)\Big|^t, \; l = 1,2 \qquad (27)$$

Applying (27) $n$ times we get (2).

Let us show that (3) is valid. Let $2 < t < 4$, $1 \le k \le n$, $U_1,\ldots,U_{k-1}$, $X_{k+1},\ldots,X_n$, $Y_1,\ldots,Y_n$ be independent symmetric r.v.'s, $E|X_i|^t < \infty$, $i = k+1,\ldots,n$, $E|Y_i|^t < \infty$, $i = 1,\ldots,n$, $a_i \ge 0$, $b_i \ge 0$, $a_i^t \le b_i$, $i = 1,\ldots,k$. Denote by $K_1$ the set of symmetric r.v.'s $X_k$ independent of $U_1,...,U_{k-1}$, $X_{k+1},...,X_n, Y_1,..., Y_n$ and satisfying the conditions $EX_k^2 \le a_k^2$, $E|X_k|^t \le b_k$, and by $K_2$ the subset of $K_1$ consisting of r.v.'s $X_k$ such that $EX_k^2 = a_k^2$, $E|X_k|^t = b_k$. From Lemma 6 with $c_i = 0$, $i = 1,...,n$, it follows that

$$\sup_{X_k \in K_l} \Big(\sum_{i=1}^{k-1} (b_i - a_i^t) E\Big|\sum_{\substack{j=1\\ j\ne i}}^{n} Y_j\Big|^t + E\Big|\sum_{i=1}^{k-1} a_i U_i (\sum_{\substack{j=1\\ j\ne i}}^{n} Y_j) + \sum_{i=k}^{n} X_i (\sum_{\substack{j=1\\ j\ne i}}^{n} Y_j)\Big|^t\Big) =$$

$$= \sum_{i=1}^{k} (b_i - a_i^t) E\Big|\sum_{\substack{j=1\\ j\ne i}}^{n} Y_j\Big|^t + E\Big|\sum_{i=1}^{k} a_i U_i (\sum_{\substack{j=1\\ j\ne i}}^{n} Y_j) + \sum_{i=k+1}^{n} X_i (\sum_{\substack{j=1\\ j\ne i}}^{n} Y_j)\Big|^t, \quad l = 1, 2 \qquad (28)$$

Using (28) $n$ times we obtain

$$\sup_{(X,n) \in M_k(n,a,b)} E\Big|\sum_{1 \le i \ne j \le n} X_i Y_j\Big|^t$$

$$= \sum_{i=1}^{n} (b_i - a_i^t) E\Big|\sum_{j=1, j\ne i}^{n} Y_j\Big|^t + E\Big|\sum_{1 \le i \ne j \le n} a_i U_i Y_j\Big|^t, \quad k = 1, 2 \qquad (29)$$

Let $1 \le k \le n$, $U_1,...,U_n, V_1,...,V_{k-1}, Y_{k+1},..., Y_n$ be independent symmetric r.v.'s, $E|Y_i|^t < \infty$, $i = k+1,..., n$, $a_i \ge 0$, $b_i \ge 0$, $a_i^t \le b_i$, $i = 1,...,n$, $c_i \ge 0$, $d_i \ge 0$, $c_i^t \le d_i$, $i = 1,...,k$. Denote by $B_1$ the set of symmetric r.v.'s $Y_k$ independent of $U_1,...,U_n, V_1,...,V_{k-1}, Y_{k+1},..., Y_n$ and satisfying the conditions $EY_k^2 \le c_k^2$,

$E|Y_k|^t \le d_k$, and by $B_2$ the subset of $B_1$ consisting of r.v.'s $Y_k$ such that $EY_k^2 = c_k^2$, $E|Y_k|^t = d_k$.

Applying Lemma 6 again with $c_i = b_i - a_i^t$ we obtain

$$\sup_{Y_k \in B_l} \Big( \sum_{i=1}^{n} (b_i - a_i^t)( \sum_{j=1,\ j\ne i}^{k-1} (d_j - c_j^t)) + \sum_{i=1}^{n} (b_i - a_i^t) E\Big| \sum_{j=1,\, j\ne i}^{k-1} c_j V_j + \sum_{j=k,\, j\ne i}^{n} Y_j \Big|^t +$$

$$+ E\Big| \sum_{j=1}^{k-1} c_j V_j ( \sum_{i=1,\ i\ne j}^{n} a_i U_i ) + \sum_{j=k}^{n} Y_j ( \sum_{i=1,\ i\ne j}^{n} a_i U_i ) \Big|^t \Big) =$$

$$= \sum_{i=1}^{n} (b_i - a_i^t)( \sum_{j=1,\, j\ne i}^{k} (d_j - c_j^t)) + \sum_{i=1}^{n} (b_i - a_i^t) E\Big| \sum_{j=1,\, j\ne i}^{k} c_j V_j + \sum_{j=k+1,\, j\ne i}^{n} Y_j \Big|^t +$$

$$+ E\Big| \sum_{j=1}^{k} c_j V_j ( \sum_{i=1, i\ne j}^{n} a_i U_i ) + \sum_{j=k+1}^{n} Y_j ( \sum_{i=1, i\ne j}^{n} a_i U_i ) \Big|^t, \quad l=1,2 \qquad (30)$$

Using (30) $n$ times we get (3).

Relations (8) and (9) might be proven in the same way.

**Proofs of theorems 4-8.** Let us prove (14). Let $2 < t < 4, D \ge 0$, and let $L(D)$ be a class of independent identically distributed r.v.'s $X_1,...,X_n$, for which

$$max(C_n^2 (E|X_1|^t)^2, (C_n^2)^{t/2} (EX_1^2)^t) = D.$$

It is evident that

$$\sup_{(X,n)\in M_1(n,D^{1/2t}/(C_n^2)^{1/4},D^{1/2}/(C_n^2)^{1/2})} E\left|\sum_{1\le i<j\le n} X_i X_j\right|^t$$

$$\le \sup_{(X,n)\in L(D)} E\left|\sum_{1\le i<j\le n} X_i X_j\right|^t$$

$$\le \sup_{(X,n)\in M_2(n,D^{1/2t}/(C_n^2)^{1/4},D^{1/2}/(C_n^2)^{1/2})} E\left|\sum_{1\le i<j\le n} X_i X_j\right|^t \qquad (31)$$

From relation (2) and its proof it follows that

$$\sup_{(X,n)\in M_k(n,D^{1/2t}/(C_n^2)^{1/4},D^{1/2}/(C_n^2)^{1/2})} E\left|\sum_{1\le i<j\le n} X_i X_j\right|^t =$$

$$=(C_n^2(1/(C_n^2)^{1/2}-1/(C_n^2)^{t/2})^2$$

$$+(1/(C_n^2)^{1/2}-1/(C_n^2)^{t/2})n/(C_n^2)^{t/4}E\left|\sum_{i=2}^{n} U_i\right|^t$$

$$+E\left|\sum_{1\le i<j\le n} U_i U_j/(C_n^2)^{1/2}\right|^t)D,\ \ k=1,2 \qquad (32)$$

(14) now follows from (31), (32) and the equality

$$B_4^*(t,n)=\sup_{D>0}(\sup_{(X,n)\in L(D)} E\left|\sum_{1\le i<j\le n} X_i X_j\right|^t/D)$$

The remaining relations (15)-(21) might be proven in the similar way.

**Acknowledgements.** We are grateful to an anonymous referee and Victor de la Pena for many useful suggestions, which led to the improvement of the paper. We would like also to thank Jaksa Cvitanic, Victor de la Pena and other participants of the probability seminar in the Department of Statistics at Columbia University for an opportunity to present and discuss the results.